\documentclass[10pt,leqno]{article}
\usepackage{amssymb,amsbsy,amsmath,amsfonts,amssymb,amscd, mathrsfs}

\usepackage[english]{babel}
\usepackage[T1]{fontenc}
\usepackage{indentfirst}

\makeatletter
\@addtoreset{equation}{section}
\makeatother

\newtheorem{statement}{}[section]
\newtheorem{theoreme}[statement]{Theorem}
\newtheorem{lemme}[statement]{Lemma}

\newtheorem{proposition}[statement]{Proposition}
\newtheorem{definition}[statement]{Definition}
\newtheorem{corollaire}[statement]{Corollary}

\newcommand\C{\mathbb C}
\newcommand\N{\mathbb N}

\newcommand\T{\mathbb T}
\newcommand\D{\mathbb D}

\newcommand\e{{\rm e}}

\newcommand\eps{\varepsilon}
\newcommand\ind{{\rm 1\kern-.30em I}}
\newcommand\qed{\hfill $\square$}

\let\phi=\varphi

\title{\bf The canonical injection of the Hardy-Orlicz space $H^\Psi$ into the Bergman-Orlicz space 
${\mathfrak B}^\Psi$}
\author{\it Pascal Lef\`evre, Daniel Li,\\ \it Herv\'e Queff\'elec, Luis Rodr{\'\i}guez-Piazza}

\date{\footnotesize \today}

\begin{document}

\maketitle

\noindent{\bf Abstract.} \emph{We study the canonical injection from the Hardy-Orlicz space $H^\Psi$ into 
the Bergman-Orlicz space ${\mathfrak B}^\Psi$.} 
\medskip

\medskip



\section{Introduction and notation}

\subsection{Introduction}

There are two natural Orlicz spaces of analytic functions on the unit disk $\D$ of the complex plane: the 
Hardy-Orlicz space $H^\Psi$ and the Bergman-Orlicz space ${\mathfrak B}^\Psi$. It is well-known that in the 
classical case  $\Psi (x) = x^p$, $H^p \subseteq {\mathfrak B}^p$ and the canonical injection $J_p$ from $H^p$ 
to ${\mathfrak B}^p$ is bounded, and even compact. In fact, for any Orlicz function $\Psi$, one has 
$H^\Psi \subseteq {\mathfrak B}^\Psi$ and the canonical injection 
$J_\Psi \colon H^\Psi \to {\mathfrak B}^\Psi$ is bounded, but we shall see in this paper that its compactness 
requires that $\Psi$ does not grow too fast. We actually characterize in 
Section~\ref{Compactness and weak-compactness} the compactness: $J_\Psi$ is compact if and only if 
$\lim_{x \to + \infty} \Psi (Ax)/ [\Psi (x)]^2 = 0$ for every $ A > 1$, and the weak compactness:  
$J_\Psi$ is weakly compact if and only if $\limsup_{x \to +\infty} \Psi (Ax)/ [\Psi (x)]^2 < +\infty$ for every 
$A > 1$ . We show that, if these two properties are ``often'' equivalent (this happens for example if $\Psi (x)/ x$ is 
non-decreasing for $x$ large enough), it is not always the case. We actually show a stronger result in 
Section~\ref{exemple}: there is an Orlicz function $\Psi$ such that $J_\Psi$ is weakly compact and is 
Dunford-Pettis, but such that $J_\Psi$ is not compact. 

\subsection{Notation}

An Orlicz function is a non-decreasing convex function $\Psi \colon [0, +\infty[ \to [0, +\infty[$ such that 
$\Psi (0) = 0$ and $\Psi (\infty) = \infty$. One says that the Orlicz function $\Psi$ has property 
$\Delta_2$ ($\Psi \in \Delta_2$) if $\Psi (2x) \leq C\, \Psi (x)$ for some constant $C > 0$ and $x$ large enough. 
It is equivalent to say that, for every $\beta > 1$,  $\Psi (\beta x) \leq C_\beta \Psi (x)$. It is known that if 
$\Psi \in \Delta_2$, then $\Psi (x) = O\, (x^p)$ for some $1 \leq p < +\infty$. One says 
(see \cite{Critere}, \cite{CompOrli}) that $\Psi$ satisfies the condition $\Delta^0$ if, for some $\beta > 1$, one 
has $\lim\limits_{x \to \infty} \Psi (\beta x) / \Psi (x) = +\infty$. If $\Psi \in \Delta^0$, then 
$\Psi (x) /x^p \mathop{\longrightarrow}\limits_{x \to \infty} + \infty$ for every $1 \leq p < \infty$. 
Indeed, let $1 \leq p < \infty$. For every $\beta > 1$ one can find $x_0 > 0$ such that 
$\Psi (\beta x)/ \Psi (x) \geq \beta^p$ for $x \geq x_0$; then 
$\Psi (\beta^n x_0) \geq \beta^{np} \Psi (x_0)$ for every $n \geq 1$. That implies that 
$\Psi (x) \geq C_p\, x^p$ for every $x > 0$ large enough. Since $p \geq 1$ is arbitrary, we get $x^p = o\,[\Psi (x)]$.
\par

We say that $\Psi \in \nabla_0 (1)$ if, for every $A > 1$, $\Psi (Ax)/ \Psi (x)$ is non-decreasing for $x$ large 
enough. This is equivalent to say (see \cite{CompOrli}, Proposition~4.7) that $\log \Psi (\e^x)$ is convex.
When $\Psi \in \nabla_0 (1)$, one has either $\Psi \in \Delta_2$, or $\Psi \in \Delta^0$.\par
\medskip

If $(S, {\cal S}, \mu)$ is a finite measure space, one defines the Orlicz space $L^\Psi (\mu)$ as the set of all 
(classes of) measurable functions $f \colon S \to \C$ for which there is a $C > 0$ such that 
$\int_S \Psi (|f|/ C)\,d\mu$ is finite. The norm $\|f\|_\Psi$ is the infimum of all $C > 0$ for which the above 
integral is $\leq 1$. The Morse-Transue space $M^\Psi (\mu)$ is the subspace of $f \in L^\Psi (\mu)$ for which 
$\int_S \Psi (|f|/ C)\,d\mu$ is finite for all $C > 0$; it is the closure of $L^\infty (\mu)$ in $L ^\Psi (\mu)$. 
One has $M^\Psi (\mu) = L ^\Psi (\mu)$ if and only if  $\Psi \in \Delta_2$.\par
If $\Psi (x)/x \mathop{\longrightarrow}\limits_{x \to + \infty} +\infty$, the conjugate function $\Phi$ of $\Psi$ is 
defined by $\Phi (y) = \sup_{x > 0} \big(xy - \Psi (x)\big)$. It is an Orlicz function and 
$[M^\Psi (\mu)]^\ast = L^\Phi (\mu)$, isomorphically.\par
We may note that if $\Psi (x) / x$ does not converges to infinity, we must have $\Psi (x) \leq ax$ for some $a \geq 1$ 
and $x$ large enough. Then $L^\Psi (\mu) = L^1 (\mu)$ isomorphically and then $\Phi (y) = +\infty$ for $y > a$ (giving 
$L^\Phi (\mu) = L^\infty (\mu)$ isomorphically).
\medskip

We denote by $\D$ the open unit disk of $\C$ and by $\T = \partial \D$ the unit circle. The normalized area-measure 
on $\D$ is denoted by ${\cal A}$ and the normalized Lebesgue measure on $\T$ is denoted by $m$.\par
The Hardy-Orlicz space $H^\Psi$ is defined as $\{ f \in H^1\,;\ f^\ast \in L^\Psi (m)\}$, where $f^\ast$ is 
the boundary values function of $f$, and $HM^\Psi = H^\Psi \cap M^\Psi (m)$ is the closure of $H^\infty$ in 
$H^\Psi$. The Bergman-Orlicz space ${\mathfrak B}^\Psi$ is the subspace of analytic $f \in L^\Psi ({\cal A})$, 
and ${\mathfrak B}M^\Psi = {\mathfrak B}^\Psi \cap M^\Psi ({\cal A})$ is the closure of $H^\infty$ in 
${\mathfrak B}^\Psi$. Since, for $f \in H^\Psi$, $\|f \|_{H^\Psi} = \sup_{0 < r < 1} \|f_r \|_{H^\Psi}$ 
(see \cite{CompOrli}, Proposition~3.1), where $f_r (z) = f (r z)$, one has:
\begin{displaymath}
\int_0^{2\pi} \Psi \bigg( \frac{|f (r \e^{it})|}{\| f  \|_{H^\Psi}}\bigg)\, \frac{dt}{2\pi} 
\leq \int_0^{2\pi} \Psi \bigg( \frac{|f (r \e^{it})| }{\|f _r\|_{H^\Psi}}\bigg)\, \frac{dt}{2\pi} \leq 1\,;
\end{displaymath}
hence:
\begin{displaymath}
\int_\D \Psi \bigg( \frac{|f (r \e^{it})| }{\|f  \|_{H^\Psi}}\bigg)\, d{\cal A} 
= \int_0^1 \bigg[ \int_0^{2\pi} 
\Psi \bigg( \frac{|f (r \e^{it})| }{\|f  \|_{H^\Psi}}\bigg)\, \frac{dt}{2\pi} \bigg]\, 2r\,dr 
\leq 1\,,
\end{displaymath}
so $f \in {\mathfrak B}^\Psi$ and $\| f \|_{{\mathfrak B}^\Psi} \leq \| f \|_{H^\Psi}$. It follows that 
$H^\Psi \subseteq {\mathfrak B}^\Psi$ and the canonical injection 
$J_\Psi \colon H^\Psi \to {\mathfrak B}^\Psi$ is bounded, and has norm $1$. Let us point out that the boundedness also 
follows from 
\cite{CompOrli}, Theorem~4.10, 2), since $J_\Psi$ is a Carleson embedding 
$J_\Psi \colon H^\Psi \to {\mathfrak B}^\Psi \subseteq L^\Psi ({\cal A}$).\par
This injection is not onto, since there are functions $f \in {\mathfrak B}^\Psi$ with no radial limit on a subset 
of $\T$ of positive measure (the proof is the same as in ${\mathfrak B}^p$: see \cite{Duren-Schuster}, \S~3.2, 
Lemma~2, page 81). 
Note that $J_\Psi$ is not an into-isomorphism: take $f_n (z) = z^n$, for every $n\in\N$; it is easy to see that 
$\{f_n\}_n$ tends to $0$ in ${\mathfrak B}^\Psi$, but not in $H^\Psi$.
\bigskip

\noindent{\bf Acknowledgment.} This work is partially supported by a Spanish research project 
MTM~2009-08934. Part of this paper was made during an invitation of the second-named author by the 
\emph{Departamento de An\'alisis Matem\'atico} of the \emph{Universidad de Sevilla}. It is a pleasure to thanks 
the members of this department for their warm hospitality.

\section{Compactness and weak-compactness}\label{Compactness and weak-compactness}

In order to characterize the compactness and the weak-compactness of $J_\Psi$, we introduce the following 
quantity $Q_A$, $A > 1$:
\begin{equation}
Q_A = \limsup_{x\to +\infty} \frac{ \Psi (Ax)} {[\Psi(x)]^2}\, \raise 1pt \hbox{,}
\end{equation}
which will turn out to be essential.\par
\smallskip

We are going to start with the compactness.

\begin{theoreme}\label{compact}
The canonical injection $J_\Psi \colon H^\Psi \to {\mathfrak B}^\Psi$ is compact if and only if 
\begin{equation}\label{limite nulle}
\qquad \lim_{x \to + \infty} \frac{ \Psi (Ax)} {[\Psi(x)]^2} = 0 \quad \text{for every } A > 1\,. 
\end{equation}
\end{theoreme}

\noindent{\bf Remarks.} 1) Condition \eqref{limite nulle} means that $Q_A = 0$ for every $A > 1$. It is equivalent  
to say that:
\begin{equation}
\sup_{ A > 1} Q_A < +\infty.
\end{equation}
Indeed, assume that  $M := \sup_{ A > 1} Q_A < +\infty$. Let $0 < \eps \leq 1$ and $A > 1$; we can find 
$x_A = x_A (\eps) > 0$ such that $\Psi (A x/ \eps) / [\Psi (x)]^2 \leq 2M$ for $x \geq x_A$. By convexity, one has 
$\Psi (Ax) \leq \eps\, \Psi (A x/\eps)$, and hence $\Psi (A x) / [\Psi (x)]^2 \leq 2\eps M$ for $x \geq x_A$. We 
get $Q_A = 0$.
\par\smallskip

2) It is clear that condition \eqref{limite nulle} is satisfied whenever $\Psi \in \Delta_2$, but 
$\Psi (x) = \e^{[\log (x + 1)]^2} - 1$ satisfies \eqref{limite nulle} 
without being in $\Delta_2$. However, condition \eqref{limite nulle} implies that $\Psi$ cannot grow too fast. 
More precisely, we must have  
\begin{displaymath}
\Psi (x) = o\, (\e^{x^\alpha}) \quad \text{for every } \alpha > 0\,. 
\end{displaymath}
Indeed, one has $\Psi (A t) \leq [\Psi (t)]^2$ for $t \geq t_A$, and, by iteration, 
$\Psi (A^n t_A) \leq [\Psi (t_A)]^{2^n}$ for every $n \geq 1$. For every $x > 0$ large enough, taking 
$n \geq 1$ such that $A^n t_A \leq x < A^{n + 1} t_A$, we get $\Psi (x) \leq C_1\,\e^{C_2 x ^\alpha}$, with 
$\alpha = \log 2/\log A$. Since $A > 1$ is arbitrary, $\alpha$ may be any positive number. The little-oh 
condition follows from the fact that the inequality is true for all $\alpha > 0$.
\medskip

\noindent{\bf Proof of Theorem~\ref{compact}.} By definition, ${\mathfrak B}^\Psi$ is a subspace of 
$L^\Psi (\D, {\cal A})$; hence we can see $J_\Psi$ as a Carleson embedding 
$J_\Psi \colon H^\Psi \to L^\Psi (\D, {\cal A})$. If $S (\xi, h) = \{z\in \D\,;\ |z - \xi| < h\}$, the compactness 
of $J_\Psi$ implies, by \cite{CompOrli}, Theorem~4.11, that, for every $A > 1$, every $\eps > 0$, and $h > 0$ 
small enough:
\begin{displaymath}
h^2 \leq 4\,{\cal A} [S (\xi, h)] \leq \frac{4 \eps}{\Psi [ A \Psi^{- 1} (1/h)]}\, \raise 1pt \hbox{,}
\end{displaymath}
that is, setting $x = \Psi^{-1} (1/h)$, $\Psi (Ax) \leq 4 \eps\, [\Psi (x)]^2$, and \eqref{limite nulle} is 
satisfied.\par
Conversely, one has:
\begin{displaymath}
\sup_{0 < t \leq h} \sup_{|\xi|=1} \frac{{\cal A} [S (\xi, t)]} {t} \leq \sup_{0 < t \leq h} \frac{t^2}{t} = h\,,
\end{displaymath}
which is $o \big( (1/h)/ \Psi [A \Psi^{- 1} (1/h)] \big)$ for every $A > 1$, if \eqref{limite nulle} holds; hence, 
by \cite{CompOrli}, Theorem~4.11, again, $J_\Psi$ is compact.
\qed
\bigskip

We now turn ourself to the weak compactness.

\begin{theoreme}\label{weak compactness}
The following assertions are equivalent:\par
{\rm (a)} $J_\Psi \colon H^\Psi \to {\mathfrak B}^\Psi$ is weakly compact; \par
{\rm (b)} $J_\Psi$ fixes no copy of $c_0$;\par
{\rm (c)} $J_\Psi$ fixes no copy of $\ell_\infty$;\par
{\rm (d)} $Q_A < +\infty$, for every $ A > 1$;\par
{\rm (e)} $H^\Psi \subseteq {\mathfrak B}M^\Psi$;\par
{\rm (f)} $J_\Psi$ is strictly singular.
\end{theoreme}

Recall that an operator $T \colon X \to Y$ between two Banach spaces is said to be strictly singular if there is 
no infinite-dimensional subspace $X_0$ of $X$ on which $T$ is an into-isomorphism.\par
\medskip

The proof will be somewhat long, and before beginning it, we shall remark that if $\Psi \in \Delta^0$, then condition
\begin{equation}\label{limsup finie}
\qquad Q_A < +\infty \quad \text{for every } A > 1
\end{equation}
implies condition~\eqref{limite nulle}. Indeed, if 
$\lim\limits_{x \to +\infty} \frac{\Psi (\beta x)} {\Psi (x)} = + \infty$, we get, for every $A > 1$:
\begin{displaymath}
\limsup_{x \to +\infty} \frac {\Psi (Ax)} {[\Psi (x)]^2} 
=  \limsup_{x \to +\infty} \frac {\Psi (Ax)} {\Psi (\beta A x)}\, \frac {\Psi (\beta Ax)} {[\Psi (x)]^2} 
\leq \limsup_{x \to +\infty} \frac {\Psi (Ax)} {\Psi (\beta A x)}\, Q_{\beta A} = 0\,.
\end{displaymath}
Now, if, for some $A > 1$, $\Psi (Ax)/ \Psi (x)$ is non-decreasing for $x$ large enough (in particular if 
$\Psi \in \nabla_0 (1)$), one has the dichotomy: either $\Psi \in \Delta_2$, and then $J_\Psi$ is compact; or 
$\Psi \in \Delta^0$ and hence the weak compactness of $J_\Psi$ implies, by the two above theorems, its compactness. 
Hence:

\begin{proposition}
If, for some $A > 1$, $\Psi (Ax) / \Psi (x)$ is non-decreasing, for $x$ large enough, then the weak compactness of 
$J_\Psi$ is equivalent to its compactness.
\end{proposition}

However, it is easy to construct an Orlicz function $\Psi$ which satisfies condition~\eqref{limsup finie}, but not 
condition~\eqref{limite nulle}. We do not give an axample here because we have a stronger result in 
Section~\ref{exemple}.\par
\bigskip

In order to prove Theorem~\ref{weak compactness}, we shall need several lemmas.

\begin{lemme}\label{lemme 1}
Let $\Psi$ be any Orlicz function. If we define $\Psi_1 (t) = [\Psi (t) ]^2$, $t\ge 0$, then $\Psi_1$ is an 
Orlicz function for which $H^\Psi \subseteq {\mathfrak B}^{\Psi_1}$ and the canonical injection of $H^\Psi$ into
${\mathfrak B}^{\Psi_1}$ is continuous.
\end{lemme}

\noindent{\bf Proof.} It is enough to see that $H^\Psi$ continuously embeds into $L^{\Psi_1}({\cal A})$, and for 
this we can use Theorem~4.10 in \cite{CompOrli}. Following the notation of that theorem for the measure 
$\mu = {\cal A}$, it is easy to see that, as $h\to 0^+$, $\rho_{\cal A} (h) \approx h^2$, and 
$K_{\cal A} (h) \approx h$.
Observe that, for $t > 1$, we have  $\Psi_1 [\Psi^{- 1} (t)] = t^2$, and so, for $h\in (0,1)$,
\begin{displaymath}
\frac {1/h} {\Psi_1 [\Psi^{-1}(1/h) ]} = \frac {1/h} {1/h^2} = h \succeq  K_{\cal A} (h).
\end{displaymath}
Using part 2) of Theorem~4.10  in \cite{CompOrli}, the lemma follows.
\qed

\begin{lemme}\label{lemme 2}
Let $M > \delta > 0$ and $\{f_n\}_n$ be a sequence in $H^\Psi \cap {\mathfrak B}M^\Psi$ such that:\par
{\rm (a)} $\{f_n\}_n$ tends to $0$ uniformly on compact subsets of \ $\D$;\par
{\rm (b)} $\|f_n\|_{{\mathfrak B}^\Psi} \ge \delta$, for every $n \geq 1$; \par
{\rm (c)} $\|f_n\|_{H^\Psi} \le M$, for every $n \geq 1$.\par
\smallskip 
Then there exists a subsequence $\{f_{n_k}\}_k$ such that $\sum_k |f_{n_k}(z)| < +\infty$, for every 
$z\in \D$, and for every $\alpha = (\alpha_k)_k \in \ell_\infty$, one has, writing 
$T\alpha (z) =\sum_{k = 1}^\infty \alpha_k f_{n_k} (z)$:
\begin{equation}\label{isomorph}
T\alpha \in {\mathfrak B}^\Psi \qquad \text{and } \quad
 (\delta/2) \|\alpha\|_\infty \le \|T\alpha\|_{{\mathfrak B}^\Psi} \le 2M  \|\alpha\|_\infty.
\end{equation}
\end{lemme}

\noindent{\bf Remark. } It is clear that, by \eqref{isomorph}, we are defining an operator $T$ from 
$\ell_\infty$ into ${\mathfrak B}^\Psi$ which is an isomorphism between  $\ell_\infty$ and its image. 
In particular, the subsequence $\{f_{n_k}\}_k$ is  equivalent, in ${\mathfrak B}^\Psi$, to the 
canonical basis of $c_0$.
\medskip

\noindent{\bf Proof.} First we are going to construct, inductively, a subsequence 
$\{f_{n_k}\}_k$ of $\{f_n\}$, and an increasing sequence $\{r_k\}_k$ in $(0,1)$, such that 
$\lim_{k\to \infty} r_k = 1$ and, setting
\begin{displaymath}
D_k =\{ z\in \D \, ;\  |z| \le r_k\}\,, \quad \text{for } k\ge 1, 
\end{displaymath}
and
\begin{displaymath}
 C_1 = D_1\,,\qquad   
 C_k = D_k\setminus D_{k - 1} =\{ z\in \D \,;\  r_{k - 1} < |z| \le r_k\}, \quad k\ge 2, 
\end{displaymath}
we have:
\begin{equation}\label{deux}
\qquad  |f_{n_k} (z)|\le 2^{-k}, \quad \text{for every } z\in D_{k - 1}, \text{ and every } k\ge 2\,;
\end{equation}
and%
\begin{equation}\label{trois}
\qquad \|f_{n_k} \ind_{\D \setminus C_k} \|_{L^\Psi} < \delta  2^{- k - 2}\,, \quad 
\text{for every } k\ge 1.
\end{equation}

Start the construction by taking $n_1 = 1$. It is a known fact that, for every function $f$ in the 
Morse-Transue space $M^\Psi ({\cal A})$, we have
\begin{equation}\label{quatre}
\lim_{{\cal A} (A)\to 0} \|f\, \ind_A\|_{L^\Psi} = 0.
\end{equation}
Now, using \eqref{quatre}, with $f = f_{n_1}$ and  considering sets $A$ of the form 
$A =\{z\in \D \,;\  r <|z| < 1\}$, we get $r_1\in (0,1)$ so that,  for  
$C_1 = D_1=\{ z\in \D \,;\  |z| \le r_1\}$,  we have
\begin{displaymath}
\|f_1\ind_{\D \setminus C_1}\|_{L^\Psi} < \delta  2^{- 3}\,.
\end{displaymath}
By the uniform convergence of $\{f_n\}_n$ to $0$ on $D_1$, we can find $n_2 > n_1$ such that 
\begin{displaymath}
|f_{n_2}(z)|\le 1/4, \text{ for every } z\in D_1, \qquad \text{and } \qquad
 \|f_{n_2}\ind_{D_1}\|_{L^\Psi} < \delta  2^{- 5}\,.
\end{displaymath}
Using this last inequality and \eqref{quatre} again (for $f = f_{n_2}$), we get $r_2\in (r_1, 1)$, 
$r_2 > 1 -1/2$, such that, setting $C_2 = \{ z\in \D \,;\ r_1< |z| \le r_2\}$, we have
\begin{displaymath}
\|f_{n_2} \ind_{\D \setminus C_2}\|_{L^\Psi} <\delta  2^{- 4}\,.
\end{displaymath}
Now that we have \eqref{deux} and \eqref{trois} for $k = 1$ and $k = 2$, it is clear how we must 
iterate the inductive construction. At the time of choosing $r_k \in (r_{k-1}, 1)$, we also impose the 
condition $r_k > 1 - 1/k$ in order to get $\lim_{k\to \infty} r_k = 1$.
\smallskip

Once the construction is achieved, let us see why the subsequence $\{f_{n_k}\}_k$ works. 
The condition \eqref{deux} and the fact that $\lim_{k\to \infty} r_k = 1$ imply that, for every 
compact set $K$ in $\D$ and $z\in \D$, there exists $l_K\in\N$ such that:
\begin{displaymath}
|f_{n_k}(z)|\le 2^{-k}, \quad \text{for every } z\in K, \text{ and every } k\ge l_K\,.
\end{displaymath}
This yields two facts. First, $\sum_k |f_{n_k}(z)| < +\infty$, for every $z\in \D$, and 
secondly: for every bounded complex sequence $\alpha = (\alpha_k)_k \in \ell_\infty$,
the series $\sum_k \alpha_k f_{n_k}$ converges uniformly on compact subsets 
of $\D$, and its sum, the function $T\alpha$, is analytic on $\D$.
\smallskip

It remains to prove the estimates in \eqref{isomorph} about the norm of $T\alpha$ in $L^\Psi({\cal A})$. 
By homogeneity, we may assume that $\|\alpha\|_\infty = 1$. Let us write 
$g_k = f_{n_k}\ind_{C_k}$ and $h_k = f_{n_k} \ind_{\D \setminus C_k}$, for every $k\ge 1$,
\begin{displaymath}
g = \sum_{k=1}^\infty \alpha_k g_k\qquad \text{and} \qquad 
h = \sum_{k=1}^\infty \alpha_k h_k \,.
\end{displaymath}
We have $T\alpha = g + h$. By \eqref{trois} and the fact that $|\alpha_k|\le 1$, 
we have that $h	\in L^\Psi ({\cal A})$ and $\|h\|_{L^\Psi}\le \delta/4$. 
\smallskip

By the condition (c) in the statement and the definition of the norm in $H^\Psi$ we have, 
for every $n$ and every $r \in (0,1)$:
\begin{equation}\label{cinq}
\frac {1} {2\pi} \int_0^{2\pi} \Psi \bigl( |f_n (r \e^{it})|/M\bigr)\,dt \le 1\,.
\end{equation}
The function $g_k$ is $0$ outside of $C_k$, and the sequence $\{C_k\}_k$ is a partition of $\D$. Therefore:
\begin{align*}
\int_{\D} \Psi (|g|/M) \, d{\cal A} 
& = \sum_{k=1}^\infty \int_{C_k} \Psi ( |g|/M) \, d{\cal A} 
=\sum_{k=1}^\infty \int_{C_k} \Psi (|\alpha_k|\, |f_{n_k}|/M) \, d{\cal A} \\ 
&\le    \sum_{k=1}^\infty \int_{C_k} \Psi (|f_{n_k}|/M) \, d{\cal A} \,.
\end{align*}
Integrating in polar coordinates, setting $r_0 = 0$, and using \eqref{cinq}, we get:
\begin{align*}
\int_{\D} \Psi (|g|/M) \, d{\cal A}  
& \le \sum_{k=1}^\infty \int_{r_{k -1}}^{r_k} 2r \, \frac{1} {2\pi} 
\int_0^{2\pi} \Psi (|f_{n_k} (r \e^{it})|/M) \, dt\, dr \\
& \le \sum_{k =1}^\infty \int_{r_{k - 1}}^{r_k} 2r \,dr  = 1\,,
\end{align*}
and therefore $\|g\|_{L^\Psi}\le M$, and $\|T\alpha\|_{L^\Psi}\le \delta/4 + M\le 2M$.
\smallskip

On the other hand, for every $k$, we have:
\begin{displaymath}
\|g\|_{L^\Psi}\ge \|g \, \ind_{C_k}\|_{L^\Psi} 
= |\alpha_k|\| f_{n_k} - h_k\|_{L^\Psi} 
\ge |\alpha_k|\, (\delta - \delta/2^{2+k}) \ge \frac{3\delta}{4} |\alpha_k| \,.
\end{displaymath} 
Taking the supremum on $k$, we get 
$\|g\|_{L^\Psi} \ge (3\delta/4)\, \|\alpha\|_\infty = 3\delta/ 4$. 
Consequently,
\begin{displaymath}
\|T\alpha\|_{L^\Psi}\ge \|g\|_{L^\Psi} - \|h\|_{L^\Psi} 
\ge (3\delta/4) - \delta/4  \ge \delta/2\,,
\end{displaymath}
and Lemma~\ref{lemme 2} is fully proved.
\qed
\bigskip

In the following lemma we isolate the proof of the implication 
${\rm (c)} \ \Longrightarrow \ {\rm (d)}$  in the statement of Theorem~\ref{weak compactness}.

\begin{lemme}\label{lemme 3}
Assume that the Orlicz function $\Psi$ is such that, for some $A > 1$, 
\begin{equation}\label{six}
\limsup_{x\to +\infty} \frac { \Psi (Ax) } {[\Psi (x) ]^2}=+\infty
\end{equation}
Then the injection $J_\Psi \colon H^\Psi \to {\mathfrak B}^\Psi$ fixes a copy of $\ell_\infty$.
\end{lemme}

\noindent{\bf Proof.} Let us take a sequence of positive numbers $\{d_n\}_n$, and a sequence $\{\xi_n\}_n$ 
in $\T$, such that the disks $\{ D (\xi_n, d_n)\}_n$ are pairwise disjoint  in $\D$. In particular, we should 
have $\lim_{n\to \infty} d_n = 0$.\par
The convexity of $\Psi$ implies the existence of some $c > 0$ such that $\Psi(x) \ge cx$ for every $x\ge 1$. 
Given a sequence $\{\beta_n\}_n$ in $(4, +\infty)$ to be fixed later, we can find, thanks to \eqref{six}, an 
increasing sequence $\{x_n\}$ satisfying:
\begin{equation}\label{sept}
x_n  > 1,\quad \Psi (x_n) > 1, \qquad \Psi (Ax_n) > \beta_n  [\Psi(x_n) ]^2, 
\quad \text{for every } n\in\N \,.
\end{equation}
Define $y_n$ as the point in the interval $(x_n, A x_n)$ such that
\begin{equation}\label{huit}
[\Psi (y_n) ]^2 = \Psi (Ax_n) \,.
\end{equation}

Put now $h_n = 1/\Psi (y_n)$ and  $r_n = 1 - h_n$. By \eqref{sept} and \eqref{huit}, we have 
$[\Psi(y_n) ]^2 > \beta_n > 4$, and therefore $h_n\in (0,1/2)$. Define
\begin{displaymath}
u_n (z) =\Bigl(	\frac {h_n} {1 - r_n \, \overline{\xi_n} z}\Bigr)^2 \,,\qquad 
\text{and}\qquad f_n (z) = y_n\, u_n (z)\,.
\end{displaymath}
It is easy to see that $\|u_n\|_\infty = 1$, and that $\|u_n\|_{H^1}\le h_n$. 
\smallskip

The first condition imposed to $\beta_n$ is  $\beta_n > 16/d_n^2$. That gives 
$[\Psi (y_n) ]^2 > 16/d_n^2$ and $h_n < d_n/4$. Let us write $D_n$ for the disk 
$D (\xi_n, d_n)$. Observe that, for $z \in \overline\D \setminus D_n$, we have 
\begin{displaymath}
|1 - r_n \, \overline{\xi_n}  z|=
|1 - r_n + r_n \, \xi_n \overline{\xi_n}  - r_n\, \overline{\xi_n} z|
\ge r_n|\xi_n - z| - h_n \ge (1/2) d_n - h_n \ge d_n/4\,,
\end{displaymath}
and therefore, since $[\Psi(x_n) ]^2\ge \Psi (x_n) \ge c\, x_n$,
\begin{displaymath}
|f_n (z)|\le y_n \Bigl(\frac{4h_n} {d_n}\Bigr)^2 =
\frac{16 y_n} {d_n^2 [\Psi(y_n) ]^2}\le
\frac{16 Ax_n } {d_n^2 \beta_n  [\Psi(x_n) ]^2}\le
\frac{16 A} {c \,d_n^2 \beta_n  } \,\cdot
\end{displaymath}
We also impose the condition $\beta_n > 16 A n^2/c d_n^2$, and so we  have:
\begin{equation}\label{neuf}
\qquad \qquad |f_n (z)|\le \frac {1} {n^2} \,,\qquad \text{for } z\in \overline\D \setminus D_n \,.
\end{equation}

From \eqref{neuf} we deduce that $\{f_n\}_n$ converges to $0$ uniformly on compact subsets of $\D$. 
Moreover \eqref{neuf} yields that, for every bounded sequence $\{\alpha_n\}_n$ of complex numbers, 
the series $\sum_{n\ge 1} \alpha_n f_n$ is uniformly convergent on compact subsets of $\D$. Let us write 
$f_n^*$ for the boundary value (on $\T = \partial \D$) of the function $f_n$. We claim that :
\begin{equation}\label{dix}
S =\sum_{n = 1}^\infty |f_n^*| \in L^\Psi (\T, m).
\end{equation}
From this, it is not difficult to deduce that, for every bounded sequence $\{\alpha_n\}_n$ of complex 
numbers, the function $\sum_{n = 1}^\infty \alpha_n f_n$ is in $H^\Psi$ and, for 
$M =\|S\|_{L^\Psi (\T)}$,
\begin{equation}\label{onze}
\Bigl\| \sum_{n = 1}^\infty \alpha_n f_n\Bigr\|_{H^\Psi} \le
M \|\{\alpha_n\}_n\|_\infty\,.
\end{equation}

On the other hand, taking $A_n = \{z\in \D \,;\ |z - \xi_n| \le h_n\}$, there exists a constant
$\gamma\in (0,1)$ such that ${\cal A} (A_n) \ge \gamma h_n^2$, and, for every $z\in A_n$, 
we have:
\begin{displaymath}
|1 - r_n\, \overline{\xi_n} z|\le  |1 - r_n| + |r_n\, \xi_n \overline{\xi_n} - r_n\, \overline{\xi_n} z|
=h _n + r_n\, |z-\xi_n|\le 2h_n\,,
\end{displaymath}
and consequently $|u_n (z)|\ge 1/4$. If $\delta =\gamma/4A$, we have, for every $n$,
\begin{align*}
\int_\D \Psi \Bigl( \frac {|f_n|} {\delta}\Bigr)\, d{\cal A} 
& \ge \int_{A_n} \Psi \Bigl( \frac{y_n } {4 \delta}\Bigr)\, d{\cal A} 
\ge \gamma h_n^2  \Psi \Bigl( \frac {1} {\gamma}\,A y_n\Bigr) \\
& \ge h_n^2  \Psi (A y_n) > h_n^2 \Psi (Ax_n) = 1\,.
\end{align*}
Thus $	\|f_n\|_{{\mathfrak B}^\Psi} \ge \delta$, for every $n\in\N$. We can apply 
Lemma~\ref{lemme 2}. Using this lemma and \eqref{onze}, we get a subsequence 
$\{f_{n_k}\}_k$ such that, for every $\alpha = (\alpha_k)_k \in \ell_\infty$, we have:
\begin{displaymath}
(\delta/2) \, \|\{\alpha_k\}_k\|_\infty 
\le \Bigl\| \sum_{k = 1}^\infty \alpha_k f_{n_k} \Bigr\|_{{\mathfrak B}^\Psi} 
\le \Bigl\| \sum_{k = 1}^\infty \alpha_k f_{n_k}\Bigr\|_{H^\Psi} 
\le M\|\{\alpha_k\}_k\|_\infty\,.
\end{displaymath}
This clearly says that $J_\Psi$ fixes a copy of $\ell_\infty$.
\smallskip

It remains to prove \eqref{dix}. For obtaining this we impose the last condition to the sequence 
$\{\beta_n\}_n$. We shall need:
\begin{equation}\label{douze}
\sum_{n = 1}^\infty 1/\sqrt{\beta_n} \le 1\,.
\end{equation}

Let us set $g_n =|f_n^*|\, \ind_{D_n}$. Thanks to \eqref{neuf},
$S - \sum_{n = 1}^\infty g_n$ is a bounded function. Thus we just need to prove that
$G = \sum_{n = 1}^\infty g_n$ is in $ L^\Psi (\T)$. We have  
$\|G\|_{ L^\Psi (\T) }\le A$. Indeed, recalling that the $D_n$'s are pairwise disjoint, and that
each $g_n$ is $0$ out of $D_n$, we have:
\begin{align*}
\int_\T  \Psi \Bigl(\frac {G} {A}\Bigr)\,dm 
& = \sum_{n = 1}^\infty \int_{D_n \cap \T}  \Psi \Bigl(\frac{G} {A}\Bigr)\,dm =
\sum_{n = 1}^\infty \int_{D_n \cap\T}  \Psi \Bigl(\frac {|f_n^*|} {A}\Bigr)\,dm \\ 
& \le \sum_{n = 1}^\infty \int_{\T}  \Psi\Bigl(\frac{y_n|u_n^*|} {A}\Bigr)\,dm \\ 
\noalign{\noindent and by the convexity of $\Psi$, and the fact that $ |u_n|\le 1$,}
&\le \sum_{n = 1}^\infty \int_{\T}  |u_n^*| \Psi \Bigl(\frac{y_n} {A}\Bigr)\,dm 
= \sum_{n = 1}^\infty \|u_n\|_{H_1} \Psi \Bigl( \frac{y_n} {A}\Bigr) \\ 
&\le \sum_{n = 1}^\infty \frac {\Psi (y_n / A)} {\Psi (y_n)} 
\le \sum_{n = 1}^\infty \frac {\Psi (x_n )} {\Psi (y_n)} 
= \sum_{n = 1}^\infty \frac {\Psi(x_n )} {\sqrt{ \Psi (A x_n)}} 
\le \sum_{n = 1}^\infty \frac { 1}{\sqrt{\beta_n}} \le 1\,, 
\end{align*}
by the required condition \eqref{douze}, and that ends the proof of Lemma~\ref{lemme 3}.
\qed
\bigskip

We are now in position to prove Theorem~\ref{weak compactness}.\par
\medskip

\noindent{\bf Proof of Theorem~\ref{weak compactness}.} We shall prove that:
\begin{displaymath}
{\rm (a)} \quad \Longrightarrow \quad 
{\rm (b)} \quad \Longrightarrow \quad 
{\rm (c)} \quad \Longrightarrow \quad 
{\rm (d)} \quad \Longrightarrow \quad 
{\rm (e)} \quad \Longrightarrow \quad 
{\rm (a)} \,,
\end{displaymath}
and that ${\rm (b)} \Longleftrightarrow {\rm (f)}$.\par
\smallskip

The implications ${\rm (a)} \Longrightarrow {\rm (b)} \Longrightarrow {\rm (c)}$ and 
${\rm (f)} \Longrightarrow {\rm (b)}$ are trivial, and we have seen in Lemma~\ref{lemme 3} 
that ${\rm (c)} \Longrightarrow {\rm (d)}$.\par
\smallskip

${\rm (d)} \Longrightarrow {\rm (e)}$. By Lemma~\ref{lemme 1}, there exists a constant $C>0$ such that, 
for every $f$ in the unit ball of $H^\Psi$, we have: 
\begin{equation}\label{treize}
\int_\D [\Psi (|f|/C)]^2 \, d{\cal A} \le 1\,.
\end{equation}
For every $A > 0$, there exist $x_A$, such that $\Psi (Ax) \le (Q_A + 1) [\Psi(x)]^2$, for every 
$x \ge x_A$. Thus for every $x \ge 0$ we have $\Psi (Ax) \le (Q_A + 1) [\Psi (x)]^2 + \Psi (A x_A)$.
Then, by \eqref{treize}, we have
\begin{displaymath}
\int_\D \Psi (A|f|/ C) \, d{\cal A} < +\infty\,,\qquad \text{for every } A>0 \,.
\end{displaymath}
Therefore $f \in {\mathfrak B}M^\Psi$, for every $f$ in the unit ball of $H^\Psi$, and thus for every 
$f$ in  $H^\Psi$.
\smallskip

${\rm (e)} \Longrightarrow {\rm (a)}$. Let $\{f_n\}_n$ be in the unit ball of $H^\Psi$. We have to 
prove that $\{f_n\}_n$ has a subsequence which converges in the weak topology of 
${\mathfrak B}^\Psi$. By Montel's Theorem $\{f_n\}_n$ has a subsequence converging uniformly on compact 
subsets of $\D$, to a function $g$ which, by Fatou's lemma, also belongs to the unit ball of $H^\Psi$. If this 
subsequence converges to $g$ in the norm of ${\mathfrak B}^\Psi$ we are done. If not, after perhaps a new 
extraction of subsequence, there exist $\delta > 0$ and a subsequence $\{f_{n_k}\}_k$, such that
\begin{displaymath}
\|f_{n_k} - g\|_{{\mathfrak B}^\Psi} \ge \delta,\qquad \text{and}\qquad 
\|f_{n_k}-g\|_{H^\Psi}\le 2\,.
\end{displaymath}
Since moreover $\{f_{n_k} - g\}_k$ converges to $0$ uniformly on compact subsets of $\D$ and,  by condition 
{\rm (e)}, $f_{n_k} - g\in {\mathfrak B}M^\Psi$, we may apply Lemma~\ref{lemme 2} and we get that  
$\{f_{n_k} - g\}_k$ has a subsequence equivalent to the canonical basis of $c_0$ in ${\mathfrak B}^\Psi$, and 
is therefore weakly null. This yields that $\{f_n\}_n$ has a subsequence converging to $g$ in the weak 
topology of  ${\mathfrak B}^\Psi$.
\smallskip

${\rm (b)} \Longrightarrow {\rm (f)}$. Suppose there exists an infinite-dimensional subspace $X$ of $H^\Psi$ 
on which the norms $\|\cdot\|_{{\mathfrak B}^\Psi}$ and $\|\cdot\|_{H^\Psi}$ are equivalent. We shall 
have finished if we prove that $X$ contains a subspace isomorphic to $c_0$ because then $J_\Psi$ will fix a 
copy of $c_0$. \par
We can assume that $X$ is contained in ${\mathfrak B}M^\Psi$ because we already know that {\rm (b)} 
implies  {\rm (e)}.
$X$ being infinite-dimensional, there exists, for every $n\in\N$,  $f_n\in X$, such that
$\|f_n\|_{H^\Psi} = 1$, and $\widehat {f_n} (k) = 0$, for $k = 0, 1, \ldots, n$. By the equivalence of the 
norms in $X$, there exists $\delta > 0$ such that $\|f_n\|_{{\mathfrak B}^\Psi} \ge \delta$, for every $n$. 
The unit ball of $H^\Psi$ is compact in the topology of ${\cal H}(\D)$. Since 
\begin{displaymath}
\lim_{n \to\infty} \widehat {f_n} (k) = 0\,,\qquad \text{for every } k\ge 0 \,,
\end{displaymath}
the only possible limit of a subsequence of $\{f_n\}_n$ is the function $0$.
So $\{f_n\}_n$ converges to $0$ uniformly on compact subsets of $\D$. As 
$f_n \in X \subseteq {\mathfrak B}M^\Psi$, for every $n$,
we can apply Lemma~\ref{lemme 2}, and we get that $\{f_n\}_n$ has a subsequence generating
an space $Y$ isomorphic to $c_0$ in ${\mathfrak B}^\Psi$. This space $Y$ is contained in 
$X$, where the norms are equivalent, so $Y$ is also isomorphic to $c_0$ for the norm of $H^\Psi$.
\qed
\par
\bigskip


\section{Other properties}

\subsection{Dunford-Pettis}

Recall that an operator $T \colon X \to Y$ between two Banach spaces $X$ and $Y$ is said to be 
\emph{Dunford-Pettis} if $\{T x_n\}_n$ converges in norm whenever $\{x_n\}_n$ converges weakly. Every 
compact operator is Dunford-Pettis. The next proposition shows that, in ``most'' of the cases, these two properties 
are equivalent for $J_\Psi$. 

\begin{proposition}\label{prop Dunford-Pettis}
If the conjugate function of $\Psi$ satisfies condition $\Delta_2$, then 
$J_\Psi \colon H^\Psi \to {\mathfrak B}^\Psi$ is Dunford-Pettis if and only if it is compact.
\end{proposition}

We shall see in Section~\ref{exemple} that without condition $\Delta_2$ for the conjugate function, $J_\psi$ 
may be Dunford-Pettis without being compact.
\medskip

\noindent{\bf Proof.} Remark first that speaking of the conjugate function of $\Psi$ implicitly assume that 
$\Psi (x)/ x$ tends to $+\infty$ as $x$ goes to $+\infty$.\par
Assume that $J_\Psi$ is not compact. By Theorem~\ref{compact}, there are some $A > 1$ and a sequence 
$\{x_j\}_j$ going to $+\infty$ such that $\Psi (A x_j) \geq [\Psi (x_j)]^2$. Setting 
$r_j = 1 - 1/ \Psi (x_j)$, this is equivalent to say that 
$A \Psi^{-1} \big(1/ (1 - r_j) \big) \geq \Psi^{-1} \big( 1/ (1 - r_j)^2\big)$. Define:
\begin{displaymath}
f_j (z) = x_j \, \bigg( \frac{1 - r_j}{1 - r_j z}\bigg)^2 \, \cdot
\end{displaymath}
One has $f_j \in HM^\Psi$ and $\| f_j\|_{H^\Psi} \leq 1$ (see \cite{CompOrli}, Corollary~3.10). Since 
$\{f_j\}_j$ converges to $0$ uniformly on compact subsets of $\D$, $\{f_j\}_j$ converges to $0$ in the 
weak-star topology of $H^\Psi$ (\cite{CompOrli}, Proposition~3.7). But, since the conjugate function of 
$\Psi$ satisfies condition $\Delta_2$, $H^\Psi$ is the bidual of $HM^\Psi$ (\cite{CompOrli}, Corollary~3.3); 
hence $\{f_j\}_j$ converges weakly to $0$ in $HM^\Psi$.\par
On the other hand, if $S_j = D (1, 1 - r_j) \cap \D$, one has $| 1 - r_j z| \leq 2 (1 - r_j)$ for $z \in S_j$; 
hence, writing $K = \| f_j\|_{{\mathfrak B}^\Psi}$, one has:
\begin{displaymath}
1 = \int_\D \Psi ( |f_j|/ K)\, d{\cal A} \geq \int_{S_j} \Psi ( |f_j|/ K)\, d{\cal A} 
\geq {\cal A} (S_j) \Psi (x_j/4 K)\,.
\end{displaymath}
Since ${\cal A} (S_j) \geq \alpha (1 - r_j)^2$, with $0 < \alpha < 1$, we get 
(since $\Psi (\alpha x_j / 4 K) \leq \alpha \Psi (x_j / 4 K)$, by convexity):
\begin{displaymath}
\| f_j\|_{{\mathfrak B}^\Psi} \geq (\alpha / 4) \, \frac{x_j}{\Psi^{-1} \big(1/ (1 - r_j)^2\big)} = 
(\alpha / 4) \, \frac{\Psi^{-1} \big( 1/(1 - r_j)\big)}{\Psi^{-1} \big(1/ (1 - r_j)^2\big)} 
\geq \frac{\alpha}{4 A}\,\cdot
\end{displaymath}
Therefore $J_\Psi$ is not Dunford-Pettis.
\qed
\bigskip

On the other hand, one has:

\begin{proposition}
If $J_\Psi$ is Dunford-Pettis, then $J_\Psi$ is weakly compact.
\end{proposition}

\noindent{\bf Proof.} By Theorem~\ref{weak compactness}, if $J_\Psi$ is not weakly compact, there is a subspace 
$X_0$ of $H^\Psi$ isomorphic to $c_0$ on which $J_\Psi$ is an into-isomorphism; hence $J_\Psi$ cannot be 
Dunford-Pettis.
\qed
\medskip

We shall see in the next section that $J_\Psi$ may be weakly compact without being Dunford-Pettis.

\subsection{Absolutely summing}

Every $p$-summing operator is weakly compact and Dunford-Pettis; so it may be expected that $J_\Psi$ is  
$p$-summing for some $p < \infty$. The next results show that this is never the case as soon as $\Psi$ 
grows faster than all the power functions.\par
\smallskip
Recall that an operator $T \colon X \to Y$ between two Banach spaces $X$ and $Y$ is called 
\emph{$(p, q)$-summing} if there is a constant $C > 0$ such that
\begin{displaymath}
\Big(\sum_{k = 1}^n \| T x_k\|^p \Big)^{1/p}\leq 
C \, \sup_{\| x^\ast \|_{X^\ast} \leq 1} \Big(\sum_{k = 1}^n | x^\ast  (x_k)|^q \Big)^{1/q} \,,
\end{displaymath}
for every finite sequence $(x_1, \ldots, x_n)$ in $X$. If $q = p$, it is said $p$-summing. Every $p$-summing 
operator is $(p, q)$-summing for $q \leq p$.
\goodbreak

\begin{theoreme}\label{theo summing}
If $J_\Psi \colon H^\Psi \to {\mathfrak B}^\Psi$ is $p$-summing, then, for every $q > p$, $\Psi (x) = O\, (x^q)$ for 
$x$ large enough. Moreover, if $p < 2$, then $J_\Psi$ is compact.
\end{theoreme}

In order to prove this, we need two lemmas.

\begin{lemme}\label{lemme avec A}
If the canonical injection $I_\Psi \colon A \to {\mathfrak B}^\Psi$ is $(p,1)$-summing, where $A = A (\D)$ is the disk 
algebra, then   $\Psi (x) = O\,(x^{2p})$ for $x$ large enough.\par
In particular, $J_r \colon H^r \to {\mathfrak B}^r$ is $(p, 1)$-summing for no $p < r/2$, and, if 
$\Psi \in \Delta^0$, then $J_\Psi$ is $(p, 1)$-summing for no $p < \infty$.
\end{lemme}

Recall that the disk algebra is the space of continuous functions on $\overline \D$ which are analytic in $\D$.\par
We refer to \cite{sommant} for a detailed study of $r$-summing Carleson embeddings 
$H^r \to L^r (\mu)$. In particular, it follows from these results that 
$J_r \colon H^r \to {\mathfrak B}^r$ is $1$-summing for $1 \leq r < 2$. On the other hand, it is 
easy to see that  $J_2 \colon H^2 \to {\mathfrak B}^2$ is not Hilbert-Schmidt ({\it i.e.} not $2$-summing): for the 
canonical orthonormal basis $\{z^n\}_n$ and $\{\sqrt{n + 1} \, z^n\}_n$ of $H^2$ and ${\mathfrak B}^2$, $J_2$ is the 
diagonal operator of multiplication by $\{1 / \sqrt{n + 1} \}_n$. It also follows from \cite{sommant} that, for 
$r \geq 2$, $J_r$ is $p$-summing for no finite $p$.
\medskip

\noindent{\bf Proof.} Assume that we do not have $\Psi (x) = O\, (x^{2p})$ for $x$ large enough. Then 
$\limsup_{x \to + \infty} \Psi (x)/ x^{2p} = + \infty$. Given any $K > 0$, take $y > 0$ 
such that $\Psi (y) / y^{2p} \geq K$ and such that $h = 1/ \sqrt{\Psi (y)} \leq 1/2$. Let $N$ be the integer part of $(1/h) + 1$. Writing $\xi_j = \e^{2\pi i j/N}$, we set: 
\begin{displaymath}
u_j (z) = \frac {h^2}{[1 - (1 - h)\, \overline{\xi_j} z]^2}\, \cdot
\end{displaymath}
We have $u_j \in A (\D)$. By \cite{CompOrli}, Lemma~5.6, one has, since $h \geq 1/N$:
\begin{displaymath}
\sum_{j=0}^{N - 1} | u_j (\e^{it})| \leq N\,h^2\, \frac {1 - (1 - h)^{2N}} {[1 - (1 - h)^2][1 - (1 - h)^N]^2} 
\leq \frac {\e^2} {(1 -\e)^2} := C\,.
\end{displaymath}
Hence:
\begin{displaymath}
\sup_{\| x^\ast \|_{A^\ast} \leq 1}  \sum_{j = 0}^{N - 1} | x^\ast (u_j)| \leq C\,.
\end{displaymath}

On the other hand, it is easy to see that $| u_j (z)| \geq 1/9$ when $|z - (1 - h)\xi_j|< h$; hence, 
if $S_j = \{z\in \D\,;\ |z - (1 - h)\xi_j|< h \}$, one has, since ${\cal A} (S_j) = h^2$:
\begin{displaymath}
1 = \int_\D \Psi \Big( \frac { |u_j (z)|} { \|u_j\|_{{\mathfrak B}^\Psi}} \Big)\,d{\cal A} (z)
\geq  \int_{S_j} \Psi \Big( \frac { 1/9} { \|u_j\|_{{\mathfrak B}^\Psi}} \Big)\,d{\cal A} 
\geq h^2 \Psi \Big( \frac { 1 /9} { \|u_j\|_{{\mathfrak B}^\Psi}} \Big) \,,
\end{displaymath}
so $\|u_j\|_{{\mathfrak B}^\Psi} \geq 1 / 9 \Psi^{-1} (1/h^2)$. Since $y = \Psi^{-1} (1/h^2)$, 
one gets:
\begin{displaymath}
\sum_{j = 0}^{N - 1} \|u_j\|_{{\mathfrak B}^\Psi}^p \geq (1 / 9)^p\frac{N}{y^p} 
\geq (1 /9)^p \bigg[ \frac{ \Psi (y)}{y^{2p}} \bigg]^{1/2} \geq \frac{K^{1/2}}{9^p}  \cdot
\end{displaymath} 
\par

This yields that the $(p, 1)$-summing norm of $I_\Psi$ should be greater than $K^{1/2p}/ 9C$, and, as $K$ is 
arbitrary, that $I_\Psi$ is not $(p, 1)$-summing.
\qed
\medskip\goodbreak

\noindent{\bf Remark.} When $I_\Psi \colon A \hookrightarrow {\mathfrak B}^\Psi$ is $p$-summing, we have 
this shorter argument.
By Pietsch's factorization theorem, this $I_\Psi$ factors through $H^p$. It follows from 
\cite{CompOrli}, Theorem~4.10, that $\alpha\,h^2 \leq \rho_{{\cal A}} (h) \leq 1 / \Psi^{-1} (A / h^{1/p})$, for 
some constants $0 < \alpha < 1$ and $A > 0$, and $h$ small enough. That means that $\Psi (x) \leq C\, x^{2p}$ 
for $x$ large enough.
\par

\begin{lemme}\label{un sommant}
If the canonical injection $I_\Psi \colon A \to {\mathfrak B}^\Psi$ is $1$-summing, then $J_\Psi$ is compact.
\end{lemme}

\noindent{\bf Proof.} The canonical injection $J_1 \colon H^1 \to {\mathfrak B}^1$ (as well as $J_\Psi$ whenever 
$\Psi \in \Delta_2$) is compact. Hence we may assume that $H^\Psi$ is not $H^1$ and hence that $\Psi (x)/x$ 
tends to $+\infty$ as $x$ tends to $+\infty$.\par
Assume that $J_\Psi$ is not compact. Then, as in the proof of Proposition~\ref{prop Dunford-Pettis}, there are 
some $A > 1$ and a sequence $\{x_k\}_k$ going to $+\infty$ such that $\Psi (A x_k) \geq [\Psi (x_k)]^2$. Setting 
$h_k = 1/ \Psi (x_k)$, we define, as in the proof of Proposition~\ref{lemme avec A}:
\begin{displaymath}
u_{k, j} (z) =  \frac{h_k^2}{[1 - (1 - h_k) \overline{\xi_{k, j}} z]^2} \, \raise 1pt \hbox{,}
\end{displaymath}
where $\xi_{k, j} = \e^{2\pi i j /N_k}$, with $N_k$ the integer part of $(1/h_k) + 1$. One has $u_ {k, j} \in A$ and 
(see the proofs of the two quoted propositions):
\begin{displaymath}
\sum_{j=0}^{N_k - 1} |u_{k, j} (\e^{it})| \leq C \qquad \text{and} \qquad  
\|u_{k, j}\|_{{\mathfrak B}^\Psi} \geq \frac{\delta \alpha}{A}\,\frac{1}{\Psi^{-1} (1/h_k)} \,\cdot 
\end{displaymath}
It follows that:
\begin{displaymath}
\sum_{j=0}^{N_k - 1} \|u_{k, j}\|_{{\mathfrak B}^\Psi} 
\geq \frac{\delta \alpha}{A} \,\frac{N_k}{\Psi^{-1} (1/h_k)} 
\geq \frac{\delta \alpha}{A} \, \frac{1/h_k}{\Psi^{-1} (1/h_k)} 
= \frac{\delta \alpha}{A} \, \frac{\Psi (x_k)}{x_k} \mathop{\longrightarrow}_{k\to \infty} + \infty.
\end{displaymath}
Hence $I_\Psi$ is not $1$-summing.
\qed
\bigskip

\noindent{\bf Proof of Theorem~\ref{theo summing}.} Since $J_\Psi \colon H^\Psi \to {\mathfrak B}^\Psi$ is 
$p$-summing and the canonical injection $I_\Psi \colon A \to {\mathfrak B}^\Psi$ factors as 
$I_\Psi \colon A \to H^\Psi \to {\mathfrak B}^\Psi$, 
this injection is $p$-summing. By Lemma~\ref{lemme avec A}, $\Psi (x) = O\, (x^{2p})$ for $x$ large enough. Hence we 
have the factorization $A \to H^{2p} \to H^\Psi \to {\mathfrak B}^\Psi$. Since the first injection is $2p$-summing and 
the last one is $p$-summing, the composition is $\max (1, p_1)$-summing, with 
$\frac{1}{p_1} = \frac{1}{2p} + \frac{1}{p}$ (see \cite{DJT}, Theorem~2.22), {\it i. e.} $p_1 = \frac{2}{3}\,p$. 
If $p_1 > 1$, we can use again Lemma~\ref{lemme avec A} with $p_1$ instead of $2p$; we get that 
$\Psi (x) = O\,(x^{2p_1})$, for $x$ large enough, and that the factorization 
$I_\Psi \colon A \to H^{2p_1} \to H^\Psi \to {\mathfrak B}^\Psi$ is $\max (1, p_2)$-summing, with 
$\frac{1}{p_2} = \frac{1}{2 p_1} + \frac{1}{p}\,$. Going on the same way, we get a 
decreasing sequence $\{p_n\}_n$ such that the canonical injection $A \to {\mathfrak B}^\Psi$ is 
$\max (1, p_n)$-summing and $\frac{1}{p_{n + 1}} = \frac{1}{2 p_n} + \frac{1}{p}\,\cdot$ 
Writing $p_n = \alpha_n p$, we get $\alpha_{n + 1} =  \frac {2 \alpha_n} {2\alpha_n + 1}$; hence 
$p_n \mathop{\longrightarrow}\limits_{n \to \infty} p/2$. In particular, $\Psi (x) = O\, (x^q)$ for every $q > p$.
\par
If $p < 2$, one has $\max (1, p_n) = 1$ for $n$ large enough, and Lemma~\ref{lemme avec A} implies that $J_\Psi$ is 
compact.
\qed
\bigskip

\noindent{\bf Remark 1.} It is not clear whether $J_\Psi$ $p$-summing, with $p \geq 2$, implies that $J_\Psi$ is 
compact. However, when $ r \geq 2$, $J_r \colon H^r \to {\mathfrak B}^r$ is $p$-summing for no $p < \infty$ 
(see \cite{sommant}).\par\medskip

\noindent{\bf Remark 2.} An operator $T \colon X \to Y$ between two Banach spaces is said to be \emph{finitely strictly 
singular} (or {\it superstrictly singular}) if for every $\eps > 0$, there is an integer $N_\eps \geq 1$ such that, 
for every subspace $X_0$ of $X$ of dimension $\geq N_\eps$, there is an $x \in X_0$ such that 
$\| T x\| \leq \eps\, \| x \|$. Every finitely strictly singular operator is strictly singular. It is not difficult to see that every compact operator is finitely strictly singular and it is shown in \cite{Plichko} (see also 
\cite{Pascal}, Corollary~2.3) that every $p$-summing operator is finitely strictly singular. We do not know when $J_\Psi$ is finitely strictly singular.

\subsection{Order boundedness}

Recall that an operator $T \colon X \to Y$ from a Banach space $X$ into a Banach lattice $Y$ is said to be 
\emph{order bounded} if there is $y\in Y_+$ such that $| T x| \leq y$ for every $x$ in the unit ball of $X$. 
Since the Bergman-Orlicz space ${\mathfrak B}^\Psi$ is a subspace of the Banach lattice 
$L^\Psi (\D, {\cal A})$, we may study the order boundedness of $J_\Psi$. Actually, we are going to see 
that the natural space for the order boundedness of $J_\Psi$ is not $L^\Psi (\D, {\cal A})$, but 
the \emph{weak Orlicz space} $L^{\Psi, \infty} (\D, {\cal A})$, the definition of which we are recalling 
below (see \cite{CompOrli}, Definition~3.16).
\begin{definition}
Let $(S, {\cal S}, \mu)$ be a measure space; the \emph{weak-$L^\Psi$ space} $L^{\Psi, \infty}$ is the 
set of the (classes of) measurable functions $f \colon S \to \C$ such that, for some constant $c > 0$, one 
has, for every $t > 0$:
\begin{displaymath}
\mu (|f| > t) \leq \frac{1}{\Psi (ct)} \,\cdot
\end{displaymath}
\end{definition}

One has $L^\Psi \subseteq L^{\Psi, \infty}$ and (\cite{CompOrli}, Proposition~3.18) the equality 
$L^\Psi = L^{\Psi, \infty}$ implies that $\Psi \in \Delta^0$. On the other hand, this equality holds 
when $\Psi$ grows sufficiently; for example, if $\Psi$ satisfies the condition $\Delta^1$: 
$x \Psi (x) \leq \Psi (\alpha x)$, for some constant $\alpha > 1$ and $x$ large enough. 
\begin{proposition}
$J_\Psi \colon H^\Psi \to {\mathfrak B}^\Psi$ is always order bounded into 
$L^{\Psi, \infty} (\D, {\cal A})$.
\end{proposition}

\noindent{\bf Proof.} Since (see \cite{CompOrli}, Lemma~3.11):
\begin{equation}\label{evaluation} 
\frac{1}{4} \Psi^{-1} \Big( \frac{1}{1 - |z|} \Big) \leq   \sup_{ \| f \|_{H^\Psi} \leq 1} | f (z) | 
\leq 4 \Psi^{-1} \Big( \frac{1}{1 - |z|} \Big)\,,
\end{equation}
one has, denoting by $S (z)$ the supremum in \eqref{evaluation}, for $t$ large enough:
\begin{displaymath}
{\cal A} (|S| > t) \leq {\cal A} \big( \{z\in \D\,;\ |z| > 1 - 1/\Psi (t/4)\} \big) \leq 
\frac{2}{\Psi (t/4)} \leq \frac{1}{\Psi (t/8)}\, \raise 1pt \hbox{,}
\end{displaymath}
and the result follows.
\qed
\medskip

Since we also have, for $t$ large enough:
\begin{displaymath}
{\cal A} (|S| > t) \geq {\cal A} \big( \{z\in \D\,;\ |z| > 1 - 1/\Psi (4t)\} \big) \geq 
\frac{1}{\Psi (4 t)} \, \raise 1 pt \hbox{,}
\end{displaymath}
we get:

\begin{corollaire}
$J_\Psi$ is order bounded into $L^\Psi (\D, {\cal A})$ if and only if $L^\Psi = L^{\Psi, \infty}$. 
This is the case if $\Psi \in \Delta^1$.
\end{corollaire}

\noindent{\bf Remark.} Contrary to the compactness, or the weak compactness, which requires 
that $\Psi$ does not grow too fast, the order boundedness of $J_\Psi$ into $L^\Psi (\D, {\cal A})$ holds when 
$\Psi$ grows fast enough. Nevertheless, for $\Psi (x) = \exp [ \big( \log (x + 1) \big)^2] - 1$, $J_\Psi$ is compact 
and order bounded into $L^\Psi (\D, {\cal A})$.
\medskip

When $J_\Psi$ is weakly compact, $J_\Psi$ maps $H^\Psi$ into 
${\mathfrak B}M^\Psi$ (Theorem~\ref{weak compactness}); hence, we may ask whether $J_\Psi$ may 
be order bounded into $M^\Psi (\D, {\cal A})$; however, we have:

\begin{proposition}
$J_\Psi$ is never order bounded into $M^\Psi (\D, {\cal A})$.
\end{proposition}

\noindent{\bf Proof.} If it were the case, we should have $S \in M^\Psi (\D, {\cal A})$, and hence 
\begin{displaymath}
\int_\D \Psi \bigg[ 4 \times \frac{1}{4} \Psi^{-1} \Big( \frac{1} {1 - |z|} \Big) \bigg]\, d{\cal A} (z) < +\infty\,,
\end{displaymath}
which is false.
\qed

\goodbreak

\section{An example}\label{exemple}

\begin{theoreme}\label{theo exemple}
There exists an Orlicz function $\Psi$ such that $J_\Psi$ is weakly compact and Dunford-Pettis, but which is not 
compact.
\end{theoreme}

Note that such an Orlicz function is very irregular: $\Psi \notin \Delta_2$, $\Psi \notin \Delta^0$, so, for 
every $A > 1$,  $\Psi (A x)/ \Psi (x)$ is not non-decreasing for $x$ large enough, and the conjugate function of 
$\Psi$ does not satisfies condition $\Delta_2$.

\bigskip

The following lemma is undoubtedly well-known, but we have found no reference, so we shall give a proof. Recall that 
a sublattice $X$ of $L^0 (\mu)$ is solid if $|f| \leq |g|$ and $g \in X$ implies $f \in X$ and $\| f \| \leq \| g\|$.

\begin{lemme}\label{lemme Rosenthal}
Let $(S, {\cal S}, \mu)$ be a measure space, and let $X$ be a solid Banach sublattice of $L^0 (\mu)$, the 
space of all measurable functions. Then, for every weakly null sequence $\{f_n\}_n$ in $X$ and 
every sequence $\{A_n\}_n$ of disjoint measurables sets, the sequence $\{f_n \ind_{A_n} \}_n$ converges weakly 
to $0$ in $X$.
\end{lemme}

\noindent{\bf Proof.} If the conclusion does not hold, there are a continuous linear functional 
$\sigma \colon X \to \C$ and some $\delta > 0$ such that, up to taking a subsequence, 
$| \sigma (f_n \ind_{A_n})| \geq \delta$. Set, for every measurable set $A \in {\cal S}$:
\begin{displaymath}
\mu_n (A) = \sigma (f_n \ind_A)\,.
\end{displaymath}
Then $\mu_n$ is a finitely additive measure with bounded variation. By Rosenthal's lemma (see 
\cite{Diestel-Uhl}, Lemma~I.4.1, page~18, or \cite{Diestel}, Chapter~VII, page~82), there is an increasing 
sequence of integers $\{n_k\}_k$ such that:
\begin{displaymath}
\Big| \mu_{n_k} \Big( \bigcup_{l \neq k} A_{n_l} \Big) \Big| \leq 
|\mu_{n_k}| \Big( \bigcup_{l \neq k} A_{n_l} \Big) \leq \delta/2\,.
\end{displaymath}
Now, if $A = \bigcup_{l \geq 1} A_{n_l}$,  $\{f_{n_k} \ind_A \}_k$ is weakly null, but:
\begin{displaymath}
|\sigma (f_{n_k} \ind_A) | 
\geq |\sigma (f_{n_k} \ind_{A_{n_k}})| - |\mu_{n_k}| \Big( \bigcup_{l \neq k} A_{n_l} \Big) 
\geq \delta - \frac{\delta}{2} = \frac{\delta}{2}\,\raise 1pt \hbox{,}
\end{displaymath}
so we get a contradiction.
\qed
\medskip

\noindent{\bf Proof of Theorem~\ref{theo exemple}.} We begin by defining a sequence $\{x_n\}_n$ of positive 
numbers in the following way: set $x_1 = 4$ and, for every $n \geq 1$, $x_{n + 1} > 2 x_n $ is the abscissa of the 
second intersection point of the parabola $y = x^2$ with the straight line containing $(x_n, x_n^2)$ and 
$(2 x_n, x_n^4)$; we have $x_{n + 1} = x_n^3 - 2x_n$ (for example, $x_2 = 56$). Define 
$\Psi \colon [0, +\infty) \to [0, +\infty)$ by $\Psi (x) = 4 x$ for $0\leq x \leq 4$, and, for $n \geq 1$:
\begin{equation}\label{def Psi}
\Psi (x_n) = x_n^2\,, \qquad \Psi (2 x_n) = x_n^4\,, \qquad 
\text{$\Psi$ affine between $x_n$ and $x_{n + 1}$}\,. 
\end{equation}
Then $\Psi$ is an Orlicz function and 
\begin{equation}\label{encadrement Psi}
\qquad x^2 \leq \Psi (x) \leq x^4 \quad \text{for} \quad x \geq 4.
\end{equation}

For this Orlicz function $\Psi$, $J_\Psi$ is not compact, since $\Psi (2x) / [\Psi (x)]^2$ does not tend to $0$. 
However, $J_\Psi$ is weakly compact, because one has the factorization 
$H^\Psi \hookrightarrow H^2 \hookrightarrow {\mathfrak B}^4 \hookrightarrow {\mathfrak B}^\Psi$ 
(by \eqref{encadrement Psi} and Lemma~\ref{lemme 1}).\par\smallskip

Assume that $J_\Psi$ is not Dunford-Pettis: there exists a weakly null sequence $\{f_n\}_n$ in the unit ball of 
$H^\Psi$ which does not converges for the norm in ${\mathfrak B}^\Psi$. Then $\{f_n\}_n$ converges uniformly 
to $0$ on the compact subsets of $\D$ (since it is weakly null) and we may assume that 
$\| f_n\|_{{\mathfrak B}^\Psi} \geq \delta$ for some $\delta > 0$.  We may also assume that 
$\| f_n\|_\infty \mathop{\longrightarrow}\limits_{n \to \infty} + \infty$ because if $\{f_n\}_n$ were uniformly 
bounded, we should have $\| f_n \|_{{\mathfrak B}^\Psi} \mathop{\longrightarrow}\limits_{n \to \infty} 0$, by 
dominated convergence.\par
We are going to show that there exist a 
subsequence $\{f_{n_k}\}_k$ and pairwise disjoint measurable sets $A_k \subseteq \T$ such that 
the sequence $\{f_{n_k} \ind_{A_k} \}_k \subseteq L^\Psi (\T, m)$ is equivalent  to the canonical basis of 
$\ell_1$, whence a contradiction with Lemma~\ref{lemme Rosenthal}.\par
\medskip

It is worth to note from now that the Poisson integral ${\cal P}$ maps boundedly $L^2 (\T)$ into $L^4 (\D)$. Indeed, 
$L^2 (\T) = H^2 \oplus \overline{H^2_0}$ and the canonical injection is bounded from $H^2$ into 
${\mathfrak B}^4$, by Lemma~\ref{lemme 1}.
\medskip

We have seen in the proof of Lemma~\ref{lemme 2} that there exist a subsequence $\{f_{n_k}\}_k$ and disjoint 
measurable annuli $C_1 = \{z\in \D\,;\ |z| \leq r_1\}$ and 
$C_k =\{z \in \D\,;\ r_{k - 1} < |z| \leq r_k \}$, $k \geq 2$, with $0 < r_1 < r_2 < \cdots < r_n < \cdots < 1$, such 
that $\| f_{n_k} \ind_{C_k} \|_{L^\Psi (\D)} \geq \delta/2$. The assumptions of that lemma are satisfied here: 
$\| f_n\|_{H^\Psi} \leq 1$, $\| f_n\|_{{\mathfrak B}^\Psi} \geq \delta$, $\{f_n\}_n$ converges uniformly 
to $0$ on the compact subsets of $\D$, and $f_n \in {\mathfrak B}M^\Psi$ because 
$H^\Psi \subseteq {\mathfrak B}M^\Psi$, since $J_\Psi$ is weakly compact. Then:\par
\medskip

\noindent{\bf Fact 1.} \emph{There exist two sequences $\{\alpha_k\}_k$ and $\{\beta_k\}_k$, with 
$\beta_n >  \alpha_n \mathop{\longrightarrow}\limits_{n \to \infty} +\infty$ such that, if 
$g_k = f_{n_k}^\ast \ind_{\{\alpha_k \leq |f_{n_k}^\ast| \leq \beta_k\}}$, then:
\begin{displaymath}
\| {\cal P} (g_k) \|_{L^\Psi (\D)} \geq \delta/3\,,
\end{displaymath}
where $f_{n_k}^\ast$ is the boundary value of $f_{n_k}$ on $\T$.}\par
\medskip

\noindent{\bf Proof.} 1) Let $\alpha_k = \frac{\delta}{12} \, \Psi^{-1} \big(1 / {\cal A} (C_k) \big)$ and 
$v_k = {\cal P} \big( f_{n_k}^\ast \, \ind_{ \{ |f_{n_k}^\ast| < \alpha_k\} } \big) \, \ind_{C_k}$. One has:
\begin{displaymath}
\int_\D \Psi \big( |v_k| / (\delta/12) \big) \, d{\cal A} 
= \int_{C_k} \Psi \big( |v_k| / (\delta/12) \big) \, d{\cal A} 
\leq \Psi \big( \alpha_k / (\delta / 12) \big) \, {\cal A} (C_k) = 1\,,
\end{displaymath}
so $\| v_k \|_{L^\Psi (\D)} \leq \delta /12$. Since 
${\cal P} (f_{n_k}^\ast ) = f_{n_k}$, we have 
$\| {\cal P} (f_{n_k}^\ast )\, \ind_{C_k} \|_{L^\Psi (\D)} 
= \| f_{n_k} \, \ind_{C_k} \|_{L^\Psi (\D)} \geq \delta / 2$, 
and we get:
\begin{displaymath}
\| {\cal P} (f_{n_k}^\ast \ind_{ \{|f_{n_k}^\ast \geq \alpha_k \} } ) \, \ind_{C_k} \|_{L^\Psi (\D)} 
\geq  \| f_{n_k} \, \ind_{C_k} \|_{L^\Psi (\D)} 
- \| v_k \|_{L^\Psi (\D)} 
\geq \frac {\delta}{2} - \frac {\delta}{12} = \frac {5\delta}{12} \,\cdot
\end{displaymath}
\par

2) Let $w_k = f_{n_k}^\ast \ind_{ \{|f_{n_k}^\ast| \geq \alpha_k \} }$. Since 
${\cal P} (w_k \, \ind_{\{|w_k| > \beta \} })$ 
tends to $0$ uniformly on $C_k$ when $\beta$ goes to infinity, Lebesgue's dominated convergence theorem gives:
\begin{displaymath}
\| {\cal P} (w_k \, \ind_{ \{|w_k| > \beta \} } ) \, \ind_{C_k}  \|_{L^\Psi (\D)} 
\leq \| {\cal P} (w_k \, \ind_{ \{|w_k| > \beta \} } ) \, \ind_{C_k} \|_{L^4 (\D)} 
\mathop{\longrightarrow}_{\beta \to +\infty} 0\,,
\end{displaymath}
so there is some $\beta_k > \alpha_k$ such that 
$\| {\cal P} (w_k \, \ind_{ \{|w_k| > \beta \} } ) \, \ind_{C_k} \|_{L^\Psi (\D)} \leq \delta / 12$.\par
We then have, with $g_k = f_{n_k}^\ast \, \ind_{ \{ \alpha_k \leq |f_{n_k}^\ast| \leq \beta_k\} }$:
\begin{displaymath}
\| {\cal P} (g_k) \|_{L^\Psi (\D)} \geq \| {\cal P} (g_k) \, \ind_{C_k} \|_{L^\Psi (\D)} 
\geq \frac{5 \delta}{12} - \frac{\delta}{12} = \frac{\delta}{3} \, \raise 1 pt \hbox{,}
\end{displaymath}
and that ends the proof of Fact~1.
\qed
\bigskip\goodbreak

\noindent{\bf Fact 2.} \emph{There are a further subsequence, denoted yet by $\{f_{n_k}\}_k$, and pairwise 
disjoint measurable subsets $E_k \subseteq \{ \alpha_k \leq |f_{n_k}^\ast| \leq \beta_k \}$, such 
that, if $h_k = f_{n_k}^\ast \,\ind_{E_k}$, then:
\begin{displaymath}
\| {\cal P} (h_k) \|_{L^\Psi (\D)} \geq \delta/4 \,.
\end{displaymath}
}

\noindent{\bf Proof.} First, since $g_k \in L^\infty (\T) \subseteq M^\Psi (\T)$, there exists $\eps_k > 0$ such that 
$m (A) \leq \eps_k$ implies $\| g_k \, \ind_A \|_{L^\Psi (\T)} \leq \delta / (12 \, \|{\cal P}\|)$ (where 
$\|{\cal P}\|$ stands for the norm of ${\cal P} \colon L^2 (\T) \to L^4 (\D)$). Now, 
${\cal P} \colon L^\Psi (\T) \to L^\Psi (\D)$ is bounded and its norm is $\leq \| {\cal P}\|$, thanks to the 
factorization $L^\Psi (\T) \hookrightarrow L^2 (\T) \hookrightarrow L^4 (\D) \hookrightarrow L^\Psi (\D)$. Hence  
$\| {\cal P} (g_k\, \ind_A) \|_{L^\Psi (\D)} \leq \delta / 12$ for $m (A) \leq \eps_k$.\par
Let $B_k = \{ \alpha_k \leq |f_{n_k}^\ast | \leq \beta_k\}$. Up to taking a subsequence, we may assume that 
$\sum_{l > k } m (B_l) \leq \eps_k$. Let 
\begin{displaymath}
E_k = B_k \setminus \bigcup_{l > k} B_l \,.
\end{displaymath}
The sets $E_k$, $k \geq 1$, are pairwise disjoint, and 
\begin{displaymath}
\| {\cal P} (g_k \, \ind_{E_k} ) \|_{L^\Psi (\D)} 
\geq \| {\cal P} (g_k \, \ind_{B_k} ) \|_{L^\Psi (\D)} 
- \| {\cal P} \big(g_k \, \ind_{\bigcup_{l > k} B_l} \big) \|_{L^\Psi (\D)} \geq \frac{\delta}{3} - \frac{\delta}{12} 
= \frac{\delta}{4} \,; 
\end{displaymath}
so we get the Fact~2 with $h_k = g_k\, \ind_{E_k} = f_{n_k}^\ast \, \ind_{E_k}$\,.
\qed

\bigskip

Set
\begin{displaymath}
F_k = \{ z \in E_k \,;\ \Psi \big( |f_{n_k}^\ast (z)\big)| \leq M\, |f_{n_k}^\ast (z) |^2 \}\,.
\end{displaymath}

For $z \in E_k \setminus F_k$, one has:
\begin{displaymath}
\int_{E_k \setminus F_k} | f_{n_k}^\ast|^2\, dm  
\leq \frac{1}{M} \int_\T \Psi ( |f_{n_k}^\ast)| \, dm  \leq \frac{1}{M} \, \raise 1 pt \hbox{,}
\end{displaymath}
so $\| f_{n_k}^\ast \,\ind_{E_k \setminus F_k} \|_{L^2 (\T)} \leq 1/ \sqrt M$ and:
\begin{align*}
\| {\cal P} (f_{n_k}^\ast \,\ind_{E_k \setminus F_k}) \|_{L^\Psi (\D)} 
& \leq \| {\cal P} (f_{n_k}^\ast \,\ind_{E_k \setminus F_k}) \|_{L^4 (\D)} \\
& \leq \| {\cal P} \| \, \|(f_{n_k}^\ast \,\ind_{E_k \setminus F_k}) \|_{L^2 (\T)} 
\leq \frac{\| {\cal P} \|}{\sqrt M} \leq \frac{\delta}{8} \, \raise 1 pt \hbox{,}
\end{align*}
for $M$ large enough. It follows that, for $M$ large enough, 
$\| {\cal P} (f_{n_k}^\ast \,\ind_{F_k}) \|_{L^\Psi (\D)} \geq \delta / 8$ and
\begin{equation}\label{mino norme}
\| f_{n_k}^\ast \, \ind_{F_k} \|_{L^\Psi (\D)} \geq  \delta/ (8\,\|{\cal P} \| ) \,.
\end{equation}
Now, we may assume that, for some $\alpha > 0$, 
\begin{displaymath}
\int_\T |f_{n_k}^\ast|^2 \, \ind_{F_k} \,dm \geq \alpha\,,
\end{displaymath}
because, if not, there would be a subsequence $\{ f_{n_{k_j}}^\ast \ind_{F_{k_j}} \}_j$ converging to $0$ in 
$L^2 (\T)$; but then $\{ {\cal P} \big( f_{n_{k_j}} \ind_{F_{k_j}} \big) \}_j$ would converge to $0$ in 
${\mathfrak B}^4$, and hence in ${\mathfrak B}^\Psi$, contrary to \eqref{mino norme}. It follows, using 
\eqref{encadrement Psi}, that:
\begin{equation}\label{mino int Psi}
\int_{F_k} \Psi (|f_{n_k}^\ast|)\, dm \geq \alpha\,.
\end{equation}

The following lemma is now the key of the proof.

\begin{lemme}\label{lemme cle}
Let $\delta_n = 2 x_{n - 1}/ x_n = 2/ (x_{n - 1}^2 - 2)$. If $\Psi (x) \leq M x^2$ and $x \geq x_n$, then, for $n$ 
large enough ($n \geq N$), one has  
$\Psi (\eps x) \geq C_M\, \eps\, \Psi (x)$ for $\delta_n \leq \eps \leq 1$.
\end{lemme}

\noindent{\bf Proof.} We may assume that $x_n \leq x < x_{n + 1}$, because if 
$x_k \leq x < x_{k + 1}$ with $k \geq n$, then $\eps \geq \delta_n$ implies $\eps \geq \delta_k$.\par

Now, remark that:
\begin{equation}\label{majo intermediaire}
\qquad\qquad \qquad \quad  \frac{\Psi (y)} {\Psi (x)} \leq 4\, \frac{y}{x} \, \raise 1pt \hbox{,} \quad \quad 
\text{for}\quad 2 x_n \leq x \leq y \leq x_{n + 1}\,.
\end{equation}
Indeed, on the one hand, 
$\frac{\Psi (y) - \Psi (x_n)} {\Psi (x) - \Psi (x_n)} = \frac{y - x_n}{x - x_n} \leq \frac{y}{x/2} = 2\, \frac{y}{x}$\,;
and, on the other hand, 
$\Psi (y) - \Psi (x_n) \geq \Psi (y) - \Psi (y/2) \geq \Psi (y) - \frac{1}{2}\, \Psi (y) = \frac{1}{2}\, \Psi (y)\,$,
so
$\frac{\Psi (y)}{\Psi (x)} \leq \frac{\Psi (y)}{\Psi (x) - \Psi (x_n)}  
\leq 2 \, \frac{\Psi (y) - \Psi (x_n)} {\Psi (x) - \Psi (x_n)} \leq 4\, \frac{y}{x}\,\cdot$
\par\medskip

We shall separate three cases:\par\smallskip

1) $\eps x \leq x_n \leq x \leq 2 x_n$. Then $\eps x \geq \eps x_n$ and hence $\Psi (\eps x) \geq \Psi (\eps x_n)$. 
But $2 x_{n  - 1} \leq \eps x_n \leq x_n$, since $\eps \geq \delta_n$; hence \eqref{majo intermediaire} implies 
that $\Psi (\eps x) \geq (\eps/4)\, \Psi (x_n) = (\eps/4)\, x_n^2$. On the other hand, one has, by 
hypothesis,  $\Psi (x) \leq M x^2 \leq  M (2 x_n)^2$, so we get $\Psi (\eps x) \geq (\eps/ 16 M) \Psi (x)$.\par
\smallskip

2) $x_n \leq \eps x \leq x \leq 2 x_n$. Then, since $1 \leq 1/\eps$:
\begin{displaymath}
\frac {\Psi (x)} {\Psi (\eps x)} \leq \frac {M x^2} {\Psi (x_n)} \leq \frac {M (2 x_n)^2} {x_n^2} 
= 4 M \leq \frac {4 M} {\eps} \,\cdot
\end{displaymath}
\par

3) For $x \geq 2 x_n$, remark that the conditions $\Psi (x) \leq M x^2$ and $x \geq 2 x_n$ imply that 
$x \geq x_n^2/\sqrt M$. Indeed, if $x \geq 2 x_n$, then $\Psi (x) \geq \Psi (2 x_n) = x_n^4$, and the 
condition $\Psi (x ) \leq M x^2$ implies $x_n^4 \leq M x^2$, {\it i.e.} $x \geq x_n^2/ \sqrt M$.\par
\par\smallskip

In this case, one has  
$\eps x \geq \eps x_n^2 / \sqrt M \geq \delta_n x_n^2 / \sqrt M = 2 (x_{n - 1}/ x_n) x_n^2 / \sqrt M 
= 2 x_{n - 1} x_n / \sqrt M \geq 2 x_n$, if $x_{n - 1} \geq \sqrt M$. Hence \eqref{majo intermediaire} gives, 
for $2 x_n \leq x < x_{n + 1}$ (since then $2 x_n \leq \eps x \leq x < x_{n + 1}$):
\begin{displaymath}
\frac {\Psi (x)} {\Psi (\eps x) } \leq 4\, \frac {x} {\eps x} = \frac {4} {\eps}\,\cdot
\end{displaymath}
That ends the proof of Lemma~\ref{lemme cle}.
\qed
\par\medskip\goodbreak

Extract now a further subsequence of $\{f_{n_k}\}$, yet denoted by $\{f_{n_k}\}$, in order that (see 
Fact 1) $\alpha_k \geq x_{N + k}$. Lemma~\ref{lemme cle} holds, with $x = \Psi (|f_{n_k}^\ast (z)|)$, $z \in F_k$, 
for every $k \geq 1$; one has (since, by definition, $\Psi (|f_{n_k}|) \leq M\,|f_{n_k}|^2$ on $F_k$):
\begin{displaymath}
\qquad \int_{F_k} \Psi (\eps\, |f_{n_k}^\ast |) \,dm \geq \eps \, C/\alpha := c\,\eps \,, 
\qquad \text{for } \delta_{N +k} \leq \eps \leq 1\,.
\end{displaymath} 

The proof of Theorem~\ref{theo exemple} reaches now its end: put $u_k = f_{n_k}^\ast \, \ind_{F_k}$, and take an 
arbitrary sequence of complex numbers such that $\sum_{k \geq 1} |\lambda_k| = 1$. Let 
$\delta_0 = \sum_{k \geq N} \delta_k$. One has $\delta_0 < 1$, because we may assume that $N$ had been taken large 
enough. One gets:
\begin{align*}
\int_\T \Psi \Big( \Big| \sum_{k \geq 1} \lambda_k u_k \Big| \Big)\, dm 
& = \sum_{k \geq 1} \int_{F_k} \Psi (| \lambda_k f_{n_k} |) \,dm \\
& \geq \sum_{ |\lambda_k| \geq \delta_{N + k} } c\, |\lambda_k| 
+ \sum_{|\lambda_k| < \delta_{N + k} } \int_{F_k} \Psi (| \lambda_k f_{n_k} |) \,dm \\
& \geq \sum_{ |\lambda_k| \geq \delta_{N + k} } c\, |\lambda_k| 
= c \, \Big( 1 - \sum_{ |\lambda_k| < \delta_{N + k} } |\lambda_k| \Big) \\
& \geq c \, \Big( 1 - \sum_{k \geq N} \delta_k \Big) = c\,  (1 - \delta_0) := c_0\,.
\end{align*}
Since $c_0 < 1$, this implies, by convexity, that 
\begin{displaymath}
\Big\| \sum_{k \geq 1} \lambda_k u_k \Big\|_{L^\Psi (\T)} \geq c_0\,.
\end{displaymath}
Hence $\{u_k\}_k$ is equivalent to the canonical basis of $\ell_1$, and that achieves the proof of 
Theorem~\ref{theo exemple}.
\qed
\par\bigskip

\noindent{\bf Remarks.} 1) It follows from Theorem~\ref{theo summing} that, for this $\Psi$, $J_\Psi$ is not 
$p$-summing for $p < 4$. By modifying the definition of $\Psi$ (taking $\Psi (x_n) = x_n^{r/2}$ and 
$\Psi (2 x_n) = x_n^r$), we get, for every $4 \leq r < \infty$, an Orlicz function $\Psi$ such that $J_\Psi$ is 
Dunford-Pettis and weakly compact, without being $p$-summing for $p < r$, and without being compact. We do not know 
whether it is possible to have $J_\Psi$ $p$-summing for no finite $p$.\par
2) Let us point out that the fact that $J_\Psi$ is Dunford-Pettis does not trivially follows from its weak compactness: 
$H^\Psi$ does not have the Dunford-Pettis property. In fact, if it were the case, the weakly compact injection 
$H^\Psi \hookrightarrow H^2$ would be Dunford-Pettis, and hence also $H^4 \hookrightarrow H^2$ (since 
$H^4 \hookrightarrow H^\Psi \hookrightarrow H^2$). But it is not the case: the sequence $\{z^n\}_n$ converges 
weakly to $0$ in $H^4$, whereas it does not converges in norm to $0$ in $H^2$.
\bigskip\goodbreak

\begin{proposition}
There is an Orlicz function $\Psi$ for which $J_\Psi$ is weakly compact, but not Dunford-Pettis.
\end{proposition}

\noindent{\bf Proof.} Let us call $\Psi_0$ the Orlicz function constructed in Theorem~\ref{theo exemple}, and let 
$\Psi (x) = \Psi_0 (x^2)$. Then, with $\beta = 2$, 
$\Psi (\beta x) = \Psi_0 (4 x^2) \geq 4 \Psi_0 (x^2) = (2 \beta) \Psi (x)$; that means that the conjugate function 
of $\Psi$ satisfies $\Delta_2$.\par
$J_\Psi$ is weakly compact (since $J_\Psi$ factors as 
$H^\Psi \hookrightarrow H^4 \hookrightarrow {\mathfrak B}^8 \hookrightarrow {\mathfrak B}^\Psi$), but is not 
compact, since $[\Psi (\sqrt{x_n})]^2 = \Psi (\sqrt{2}\, \sqrt{x_n})$. Since the conjugate function satisfies 
$\Delta_2$, $J_\Psi$ is not Dunford-Pettis, by Proposition~\ref{prop Dunford-Pettis}.
\qed


\bigskip

\vbox{\noindent{\small\it 
{\rm Pascal Lef\`evre}, Univ Lille Nord de France F-59\kern 1mm 000 LILLE, FRANCE\\
UArtois, Laboratoire de Math\'ematiques de Lens EA~2462, \\
F\'ed\'eration CNRS Nord-Pas-de-Calais FR~2956, \\
F-62\kern 1mm 300 LENS, FRANCE \\ 
pascal.lefevre@euler.univ-artois.fr 
\smallskip

\noindent
{\rm Daniel Li}, Univ Lille Nord de France F-59\kern 1mm 000 LILLE, FRANCE\\
UArtois, Laboratoire de Math\'ematiques de Lens EA~2462, \\
F\'ed\'eration CNRS Nord-Pas-de-Calais FR~2956, \\
Facult\'e des Sciences Jean Perrin,\\
Rue Jean Souvraz, S.P.\kern 1mm 18, \\
F-62\kern 1mm 300 LENS, FRANCE \\ 
daniel.li@euler.univ-artois.fr
\smallskip

\noindent
{\rm Herv\'e Queff\'elec}, Univ Lille Nord de France F-59\kern 1mm 000 LILLE, FRANCE\\
USTL, Laboratoire Paul Painlev\'e U.M.R. CNRS 8524, \\
F-59\kern 1mm 655 VILLENEUVE D'ASCQ Cedex, FRANCE \\ 
queff@math.univ-lille1.fr
\smallskip

\noindent
{\rm Luis Rodr{\'\i}guez-Piazza}, Universidad de Sevilla, \\
Facultad de Matem\'aticas, Departamento de An\'alisis Matem\'atico,\\ 
Apartado de Correos 1160,\\
41\kern 1mm 080 SEVILLA, SPAIN \\ 
piazza@us.es\par}
}

\end{document}